\documentclass [11pt,oneside,a4paper,mathscr]{amsart}

\usepackage{amscd,amsmath,amssymb,euscript}
\usepackage[frame,cmtip,curve,arrow,matrix,line,graph]{xy}

\title{Stability conditions and Kleinian singularities}
\author{Tom Bridgeland}

\date{}

\newtheorem{thm}{Theorem}[section]

\newtheorem{cor}[thm]{Corollary}
\newtheorem{prop}[thm]{Proposition}
\newtheorem{lemma}[thm]{Lemma}
\newenvironment{pf}{\paragraph{Proof}}{\qed\par\medskip}
\theoremstyle{definition}

\newtheorem{remark}[thm]{Remark}

\renewcommand{\geq}{\geqslant}

\newcommand{\K}{{{K}}}

\newcommand{\D}{\operatorname{\mathcal{D}}}
\newcommand{\Coh}{\operatorname{Coh}}
\newcommand{\Aut}{\operatorname{Aut}}
\newcommand{\isom}{\cong}

\newcommand{\tensor}{\otimes}

\newcommand{\C}{\mathbb C}

\newcommand{\Z}{\mathbb Z}

\newcommand{\OO}{\mathcal O}

\newcommand{\Ext}{\operatorname{Ext}}
\newcommand{\Hom}{\operatorname{Hom}}
\newcommand{\eu}{\operatorname{\chi}}

\newcommand{\lra}{\longrightarrow}
\newcommand{\R}{\mathbb{R}}

\renewcommand{\AA}{\mathcal{A}}
\newcommand{\h}{\mathfrak{h}}
\newcommand{\g}{\mathfrak{g}}
\newcommand{\gh}{\hat{\g}}
\newcommand{\hreg}{\h^{\text{reg}}}
\newcommand{\Stab}{\operatorname{Stab}}

\newcommand{\Gammah}{\hat{\Gamma}}
\newcommand{\Dh}{\hat{\D}}
\newcommand{\AAh}{\hat{\AA}}
\newcommand{\Ah}{\hat{A}}
\newcommand{\hh}{\hat{\h}}
\newcommand{\hregh}{\hat{\h}^{\text{reg}}}
\newcommand{\Lambdah}{\hat{\Lambda}}
\newcommand{\Wh}{\hat{W}}

\renewcommand{\Im}{\operatorname{Im}}
\renewcommand{\Re}{\operatorname{Re}}
\renewcommand{\P}{\mathcal{P}}
\newcommand{\ZZ}{\mathcal{Z}}
\newcommand{\CC}{\mathcal{C}}
\newcommand{\Br}{\operatorname{Br}}
\newcommand{\SL}{\operatorname{SL}}

\begin{document}

\begin{abstract}
We describe the spaces of stability conditions on certain triangulated categories associated to Dynkin diagrams. These categories can be defined  either algebraically via module categories of preprojective algebras, or geometrically via coherent sheaves on resolutions of Kleinian singularities. The resulting spaces are related to regular subsets of the corresponding
complexified  Cartan algebras.
\end{abstract}

\maketitle

\section{Introduction}
In this paper we describe the spaces of stability conditions \cite{Br1} on certain triangulated categories associated to Dynkin diagrams. These categories can be defined  either algebraically via module categories of preprojective algebras, or geometrically via coherent sheaves on resolutions of Kleinian singularities.  
The resulting categories behave in almost all respects like derived categories of coherent sheaves on  K3 surfaces,
and the results developed in \cite{Br2} quickly yield Theorems 1.1 and 1.3 below. We give the details here since they provide good  examples of spaces of stability conditions and have some interesting connections with other areas of representation theory. The results were obtained independently by several other mathematicians including A. Ishii and H. Uehara.

\subsection{}
Let $G\subset \SL_2(\C)$ be a finite subgroup and let $\Coh_G(\C^2)$ denote the abelian category of $G$-equivariant coherent sheaves on $\C^2$. Let $\AAh\subset \Coh_G(\C^2)$ denote the full subcategory consisting of equivariant sheaves supported at the origin in $\C^2$, and let $\AA\subset \AAh$ be the full subcategory consisting of equivariant sheaves with no non-trivial $G$-invariant sections.
Let $\D$ and $\Dh$  be the full subcategories of $\D^b\Coh_G(\C^2)$ consisting of complexes whose cohomology sheaves lie in $\AA$ and $\AAh$ respectively.  The aim of this paper is to describe the spaces
of stability conditions on these triangulated categories.
Before describing the results in more detail we give some alternative descriptions of our categories.

The first description is more geometric. Let $X=\C^2/G$ be the Kleinian quotient singularity associated to $G$ and let $f\colon Y\to X$ be the minimal resolution of singularities. The derived McKay correspondence gives an equivalence
\[\D^b\Coh_G(\C^2) \lra \D^b \Coh(Y).\]
It is easy to check that under this equivalence the subcategory $\Dh$ corresponds to the full subcategory of $\D^b\Coh(Y)$ consisting of objects supported on the exceptional divisor $f^{-1}(0)\subset Y$, and the subcategory $\D$ corresponds to the full subcategory consisting of objects $E$ satisfying $\mathbf{R} f_*(E)=0$.
 
The second description is more algebraic. Recall that
J. McKay \cite{Mc} showed how to associate an extended Dynkin graph $\Gammah$ to our finite subgroup $G\subset \SL_2(\C)$. The vertices of $\Gammah$ are labelled by the ismomorphism classes of irreducible representations  of $G$, and the vertices corresponding to two irreducible representations $\rho_i$ and $\rho_j$ are joined by an edge precisely when $\rho_i\subset Q\tensor\rho_j$, where $Q$ is the given representation $G\subset SL_2(\C)$. The possible graphs $\Gammah$ are shown in Figure 1, with the special vertex corresponding to the trivial representation of $G$  marked with an open dot. Removing this vertex leaves a Dynkin graph $\Gamma$.

\begin{figure}
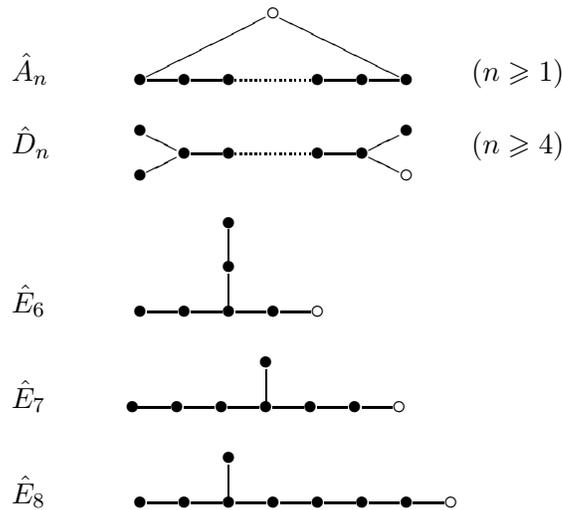

\[
\begin{array}{ll}
\hat{A}_n  \qquad &  \xy /r.07pc/:
\POS (0,0) *{\bullet} ="v1",
(60,30) *{\circ} ="v0",
(20,0) *{\bullet} ="v2",
(40,0) *{\bullet} ="v3",
(80,0) *{\bullet} ="v4",
(100,0) *{\bullet} ="v5",
(120,0) *{\bullet} ="v6"
\POS"v1" \ar@{-} "v2"
\POS"v2" \ar@{-} "v3"
\POS"v3" \ar@{.} "v4"
\POS"v4" \ar@{-} "v5"
\POS"v5" \ar@{-} "v6"
\POS"v0" \ar@{-} "v1"
\POS"v0" \ar@{-} "v6"
\endxy \qquad (n\geq 1)
\\
& \\
\hat{D}_n &  \xy /r.07pc/:
\POS (0,-10) *{\bullet} ="v1",
\POS (0,10) *{\bullet} ="v0",
(20,0) *{\bullet} ="v2",
(40,0) *{\bullet} ="v3",
(80,0) *{\bullet} ="v4",
(100,0) *{\bullet} ="v5",
(120,-10) *{\circ} ="v6",
(120,10) *{\bullet} ="v7",
\POS"v1" \ar@{-} "v2"
\POS"v0" \ar@{-} "v2"
\POS"v2" \ar@{-} "v3"
\POS"v3" \ar@{.} "v4"
\POS"v4" \ar@{-} "v5"
\POS"v5" \ar@{-} "v6"
\POS"v5" \ar@{-} "v7"
\endxy \qquad (n\geq 4) \\
& \\
\hat{E}_6 &  \xy /r.07pc/:
\POS (0,0) *{\bullet} ="v1",
(20,0) *{\bullet} ="v2",
(40,0) *{\bullet} ="v3",
(60,0) *{\bullet} ="v4",
(80,0) *{\circ} ="v5",
(40,20) *{\bullet} ="v7",
(40,40) *{\bullet} ="v8"
\POS"v1" \ar@{-} "v2"
\POS"v2" \ar@{-} "v3"
\POS"v3" \ar@{-} "v4"
\POS"v4" \ar@{-} "v5"
\POS"v3" \ar@{-} "v7"
\POS"v7" \ar@{-} "v8"
\endxy \\
& \\
\hat{E}_7 &  \xy /r.07pc/:
\POS (0,0) * ={\bullet}="v1",
(20,0) *{\bullet}="v2",
(40,0) *{\bullet}="v3",
(60,0) *{\bullet} ="v4",
(80,0) *{\bullet} ="v5",
(100,0) *{\bullet}="v6",
(120,0) *{\circ}="v7",
(60,20) *{\bullet}="v0"
\POS"v1" \ar@{-} "v2"
\POS"v2" \ar@{-} "v3"
\POS"v3" \ar@{-} "v4"
\POS"v4" \ar@{-} "v5"
\POS"v5" \ar@{-} "v6"
\POS"v6" \ar@{-} "v7"
\POS"v4" \ar@{-} "v0"
\endxy \\
& \\
\hat{E}_8 &  \xy /r.07pc/:
\POS (0,0) *{\bullet} ="v1",
(20,0) *{\bullet} ="v2",
(40,0) *{\bullet} ="v3",
(60,0) *{\bullet}="v4",
(80,0) *{\bullet}="v5",
(100,0) *{\bullet}="v6",
(120,0) *{\bullet}="v8",
(140,0) *{\circ}="v9",
(40,20) *{\bullet}="v0"
\POS"v1" \ar@{-} "v2"
\POS"v2" \ar@{-} "v3"
\POS"v3" \ar@{-} "v4"
\POS"v4" \ar@{-} "v5"
\POS"v5" \ar@{-} "v6"
\POS"v6" \ar@{-} "v8"
\POS"v3" \ar@{-} "v0"
\POS"v8" \ar@{-} "v9"
\endxy
\end{array}
\]
\caption{The affine Dynkin diagrams $\Gammah$.}
\label{tames}
\end{figure}

The category $\Coh_G(\C^2)$ is tautologically the same thing
as the category of modules for the skew group algebra $\C[x,y]*G$. In turn it is known \cite{CH} that this algebra is Morita
equivalent to the preprojective algebra $\Ah$ of the graph $\Gammah$. More precisely, to define the preprojective algebra one must choose an orientation of $\Gammah$, but different choices of orientation lead to isomorphic preprojective algebras.
Using these identifications it is easy to see that  $\AAh$ is equivalent to the category of nilpotent modules for $\Ah$, and that $\AA$ is equivalent to the full subcategory consisting of representations $M$ satisfying $e_0 M=0$, where $e_0\in \Ah$ is the idempotent corresponding to the special vertex $0$ of the quiver $\Gammah$.
From this description it is also immediate that $\AA$ is equivalent to the category of finite-dimensional modules for the preprojective algebra $A$ of the Dynkin quiver $\Gamma$.

The category $\AAh$ is finite length and has $n+1$ simple objects $S_0,\cdots,S_n$ corresponding to the vertices
of $\Gammah$. In terms of equivariant coherent sheaves these simples are of the form $S_i=\rho_i\tensor\OO_0$, where $\rho_i$ is an irreducible representation of $G$ and $\OO_0$ is the skyscraper sheaf at the origin in $\C^2$.  We shall always assume that $S_0$ corresponds to the trivial representation of $G$.
The full subcategory $\AA\subset\AAh$ consists of those objects none of whose simple factors are isomorphic to $S_0$. Clearly this is also a finite length category with simple objects $S_1,\cdots,S_n$ corresponding to the vertices of the graph $\Gamma$.

It is a slightly subtle question as to whether the category $\Dh$ is actually equivalent to $\D^b(\AAh)$ and in any case this will not matter to us. However it is worth making the point that $\D$ is definitely not equivalent to $\D^b(\AA)$. Indeed, the fact that $\Ah$ has finite global dimension implies that the category $\Dh$ is of finite type, meaning that for any objects $E$ and $F$ 
\[\dim_{\C} \big(\bigoplus_{i\in \Z} \Hom_{\Dh}(E,F[i])\big) < \infty\]
Since $\D\subset \Dh$ is a subcategory the same is true of $\D$. But the algebra $A$ has infinite global dimension, so the category $\D^b(\AA)$ is not of finite type. One could perhaps think of the category  $\D$ as a better-behaved substitute for $\D^b(\AA)$.

\subsection{}
The combinatorics of the category $\D$ are described by the root system of the finite-dimensional complex simple Lie algebra $\g$ corresponding to the Dynkin graph $\Gamma$. 
Let $\h\subset \g$ denote the Cartan subalgebra
and let $\hreg\subset \h$ be the complement of the root hyperplanes in $\h$
\[\hreg=\{\theta\in \h:\theta(\alpha)\neq 0\text{ for all }\alpha\in \Lambda\}.\]
The Weyl group $W$ is generated by reflections in the root hyperplanes and acts freely on $\hreg$.

The simple objects $S_i\in \AA$ are spherical objects in $\D$ and hence by results of P. Seidel and R.P. Thomas \cite{ST} define autoequivalences $\Phi_{S_i}\in \Aut(\D)$. We write $\Br(\D)$ for the subgroup of $\Aut(\D)$ they generate.

The following result generalises a result of Thomas \cite{Th} who proved the $A_n$ case using different methods. In fact Thomas worked with a triangulated category whose objects were dg-modules over a dg-quiver, but the formality result of \cite{ST} shows that his category is equivalent to ours (see \cite[Section 3]{Th}).

\begin{thm}
There is a connected component of the space of stability conditions $\Stab(\D)$ which is a covering space of $\hreg/W$. The subgroup $\Br(\D)\subset\Aut(\D)$  preserves this component and acts\footnote{faithfully: see Remark \ref{eend}} as the group of deck transformations.
\end{thm}

The fundamental group  of the quotient $\hreg/W$ coincides \cite{Bi,De} with the braid group $\Br(\Gamma)$ of the graph $\Gamma$. This is the group generated by elements $\sigma_1,\cdots,\sigma_n$ indexed by the vertices subject to relations
$\sigma_i\sigma_j\sigma_i=\sigma_j\sigma_i\sigma_j$
if the vertices $i$ and $j$ are connected by an edge, and $\sigma_i\sigma_j=\sigma_j\sigma_i$ otherwise.
It follows from Theorem 1.1 that there is a surjective homomorphism
\[\rho\colon \Br(\Gamma)\lra \Br(\D)\]
It follows easily from our description of $\Stab_0(\D)$ that $\rho$ sends the generator $\sigma_i$ to the twist functor $\Phi_{S_i}$.
In the $A_n$ case Seidel and Thomas \cite{ST} were able to show that $\rho$ is an isomorphism, so that $\Stab_0(\D)$ is actually the universal cover of $\hreg/W$. Reversing this argument, one might hope to find a general proof that $\Stab_0(\D)$ is simply-connected; this would then imply that $\rho$ is always an isomorphism.

Theorem 1.1 allows us to say something about the group of autoequivalences of the category $\D$. Unfortunately, we cannot rule out the possibility of exotic autoequivalences which permute the connected components of the space of stability conditions. So define $\Aut_0(\D)$ to be the subgroup of autoequivalences which preserve the connected component $\Stab_0(\D)$ of Theorem 1.1. 
A further problem is that $\Aut_0(\D)$ could in theory contain autoequivalences $\Phi$ which act trivially, in the sense that $\Phi(E)\isom E$ for all objects $E$. In fact, we are only really interested in the action of $\Aut_0(\D)$ on $\Stab_0(\D)$. So define $\Aut_0^*(\D)$ to be the quotient $\Aut_0(\D)/H$ where $H$ is the subgroup of $\Aut_0(\D)$ consisting of autoequivalences which fix every point of $\Stab_0(\D)$.

\begin{cor}Let $\Aut(\Gamma)$ be the group of symmetries of the graph $\Gamma$.
There is an isomorphism
\[\Aut_0^*(\D)\isom \Br(\D)\rtimes\Aut(\Gamma),\]
where $\Aut(\Gamma)$ acts on $\Br(\D)$ by permuting the generators $\Phi_{S_i}$. 
\end{cor}

One might wonder where the shift functor has gone in Corollary 1.2, but one can easily check by direct computation that, for example in the $A_n$ case,
one has
\[\big(\Phi_{S_1}\circ\Phi_{S_2}\circ\cdots\circ\Phi_{S_n}\big)^{n+1}=[-2],\]
with the obvious numbering of the vertices of $\Gamma$.
Presumably something similar happens in the general case.

\subsection{}
The combinatorics of the category $\Dh$ are described by the root system of the affine Kac-Moody Lie algebra $\gh$ corresponding to the graph $\Gammah$. The affine roots $\Lambdah\subset \gh^*$ span a subspace $\hh^*\subset \hat{g}^*$ which can be identified with $\h^*\oplus \C$. These roots $\Lambdah\subset \hh^*$ are of two types; the elements $(\alpha,d)$ for $\alpha\in \Lambda$ and $d\in \Z$ are the real roots; the elements $(0,d)$ for $d\in \Z\setminus\{0\}$ are the imaginary roots.

Let $\hh=\h\oplus\C$ be the dual of $\hh^*$ and let $\hregh\subset \hh$ be the complement of the affine root hyperplanes in $\hh$
\[\hregh=\{(\theta,n)\in \hh:\theta(\alpha)+nd\neq 0\text{ for all }(\alpha,d)\in \Lambdah\}.\]
It is easy to see that this is an open subset of $\hh$.
The affine Weyl group $\Wh$ is generated by reflections in the real root hyperplanes and acts freely on $\hregh$ preserving the projection to $\C^*$.

Once again, the simple objects $S_i\in \AAh$ are spherical objects in $\Dh$ and hence define autoequivalences $\Phi_{S_i}\in \Aut(\Dh)$. We write $\Br(\Dh)$ for the subgroup of $\Aut(\Dh)$ they generate. This time $\Br(\Dh)$ does not contain any power of the shift functor so we also consider the subgroup  $\Br(\Dh)\times\Z\subset\Aut(\Dh)$ where the second factor is generated by the shift functor $[2]$.

\begin{thm}
There is a connected component of the space of stability conditions $\Stab(\Dh)$ which is a covering space of $\hregh/\Wh$. The subgroup $\Br(\Dh)\times\Z\subset \Aut(\Dh)$  preserves this component and acts\footnote{faithfully: see Remark \ref{eend}} as the group of deck transformations.
\end{thm}

The fundamental group of the quotient $\hregh/\Wh$ coincides \cite{N} with the group $\Br(\Gammah)\times\Z$, where $\Br(\Gammah)$ is the braid group associated to the graph $\Gammah$. The factor $\Z$ is generated by a loop $\gamma$ around the hyperplane corresponding to the imaginary root $(0,1)$. Theorem 1.3 implies that there is a surjective homomorphism
\[\rho\colon\Br(\Gammah)\times\Z\lra \Br(\Dh)\times\Z.\]
Again it is easy to see that $\rho$ sends the generators $\sigma_i$ to the twist functors $\Phi_{S_i}$ and the generator $\gamma$ to the shift functor $[2]$.

As before define $\Aut_0(\Dh)$ to be the subgroup of autoequivalences which preserve the connected component $\Stab_0(\Dh)$  of Theorem 1.3, and $\Aut_0^*(\Dh)$ to be the quotient by the autoequivalences which act trivially on $\Stab_0(\Dh)$.  Then  we have

\begin{cor}Let $\Aut(\Gammah)$ be the group of symmetries of the graph $\Gammah$.
There is an isomorphism 
\[\Aut_0^*(\Dh)=\Z\times (\Br(\Dh)\rtimes\Aut(\Gammah)),\]
where the factor of $\Z$ is generated by the shift $[1]$, and $\Aut(\Gammah)$ acts on $\Br(\Dh)$ by permuting the generators $\Phi_{S_i}$. 
\end{cor}

A more careful analysis of the group $\Aut(\Dh)$ has been carried out in the $\hat{A}_n$ case by Ishii and Uehara \cite{IU}.

\subsection*{Acknowledgements} Thanks to Richard Thomas whose laziness \cite{Th} in only considering the $A_n$ case  made this paper possible. Thanks also to Bill Crawley-Boevey and Alastair King for useful comments, and to Yukinobu Toda for pointing out several mistakes in an earlier version.

% *******************************************************************************************************************
% *******************************************************************************************************************
% *******************************************************************************************************************

\section{Background}

\subsection{}
Here we state some simple facts about the categories defined in the introduction which will be used in the proofs of our main Theorems. These results are all well-known and we only sketch the proofs. Further details can be found in \cite{CH}.

\begin{lemma}
The categories $\D$ and $\Dh$ are triangulated categories of finite type with Serre functor $[2]$.
\end{lemma}

\begin{pf}
As explained in the introduction, we can identify $\Dh$ with the full subcategory of the derived category of coherent sheaves on the minimal resolution $f\colon Y\to\C^2/G$ consisting of objects supported on the exceptional fibre $f^{-1}(0)$. It follows that $\Dh$ has finite type. Since $Y\to \C^2/G$ is crepant the canonical bundle $\omega_Y$ is trivial on this fibre. The result then follows from Serre duality.
\end{pf}

Recall that an object $S\in \Dh$ is spherical if 
\[
\Hom^k_{\Dh}(S,S)=\biggl\{\begin{array}{ll} \C & \text{ if }k=0 \text{ or }2,
\\ 0 &\text{ otherwise.}\end{array}\biggr.
\]
It follows from constructions given in \cite{ST}
that any such object defines an
auto--equivalence $\Phi_S\in \Aut \Dh$ called a twist functor such that for any $E\in
\Dh$ there is a triangle
\[\Hom_{\Dh}(S,E)\tensor S\lra E\lra \Phi_S(E).\]
Note that at the level of the Grothendieck group the functor induces a reflection
\[\phi_S([E])=[E]-\eu([S],[E])[S].\]
These twist functors will be very important in what follows.

\begin{lemma}
\label{one}
The abelian category $\AA$ is of finite type with simple objects $S_1,\cdots, S_{n}$ labelled by the vertices of the graph $\Gamma$. Each of these objects is spherical in $\D$. Given any two of these simples the space $\Hom^1_{\D}(S_i,S_j)$ is one or zero-dimensional depending on whether the corresponding vertices of $\Gamma$ are joined by an edge or not. 
\end{lemma}

\begin{pf}
This is an easy computation of $\Ext$-groups, either in the category $\Coh_G(\C^2)$ or the category of representations of the preprojective algebra.
\end{pf}

It follows from this Lemma that $\eu(S_i,S_i)=2$ for all $i$ and $\eu(S_i,S_j)=0$ or 1 depending on whether $i$ and $j$ are connected by an edge in $\Gamma$. Thus the Grothendieck group $K(\D)$ with its Euler form can be identified with the root lattice $\Z \Lambda\subset \h^*$ of the corresponding finite-dimensional Lie algebra $\g$, equipped with the unique multiple of the Killing form for which $\alpha^2=2$ for all roots $\alpha\in \Lambda$.

Under this identification the classes of spherical objects of $\D$ correspond to the roots $\Lambda$, and the classes $[S_i]$ form a system of simple roots. The reflection $\phi_S$ of $K(\D)$ defined by a spherical object $S$ of $\D$ induces the root reflection of the root lattice defined by the corresponding element of $\Lambda$.

\begin{lemma}
\label{two}
The abelian category $\AAh$ is of finite type with simple objects $S_0,\cdots, S_{n}$ labelled by the vertices of the graph $\Gammah$. Each of these objects is spherical in $\Dh$. Given any two of these simples the space $\Hom^1_{\Dh}(S_i,S_j)$ is one or zero-dimensional depending on whether the corresponding vertices of $\Gammah$ are joined by an edge or not.
\end{lemma}

\begin{pf}
This is the same as the proof of Lemma \ref{one}.
\end{pf}

Let $\gh$ be the affine Kac-Moody Lie algebra corresponding to the graph $\Gammah$ (see \cite{Ka} for definitions). Let $\Z\Lambdah\subset \gh^*$ be the affine root lattice and let $\hh^*$ be the vector subspace of $\gh^*$ spanned by $\Lambdah$. Note that this is not the dual of the Cartan subalgebra of $\gh$.
The normalised invariant form on $\gh$ induces a form on $\hh^*$ with a one-dimensional kernel.
We can identify $\hh^*$  with the direct sum $\h^*\oplus \C$. The restriction of the invariant form can then be written $(\alpha,d)^2=\alpha^2$.

As before, Lemma \ref{two} is enough to calculate the Euler form on the Grothendieck group $K(\Dh)$. This time the form is indefinite with a one-dimensional kernel generated by the class of the equivariant sheaf $\C[G]\tensor\OO_0$. The group $K(\Dh)$ with the Euler form can be identified with the root lattice $\Z\Lambdah$ with the restriction of the invariant form. The vector space $\hh^*$ becomes identified with $K(\Dh)\tensor\C$.

Under these identifications the classes of the spherical objects of $\Dh$ correspond to the real roots $(\alpha,d)$ for $\alpha\in \Lambda\subset \h^*$ and $n\in \Z$. The class of the equivariant sheaf $\C[G]\tensor\OO_0$ corresponds to the imaginary root $(0,1)$. Once again, the reflection $\phi_S$ of $K(\Dh)$ defined by a spherical object $S$ of $\Dh$ induces the reflection of the root lattice defined by the corresponding real root.

\subsection{}
We refer to the reader to \cite{Br1} for definitions concerning stability conditions. Here we just give a brief summary,
mainly in order to fix notation.

A stability condition $\sigma=(Z,\P)$
on a triangulated category $\D$ is defined by full abelian subcategories $\P(\phi)\subset \D$
for each $\phi\in\R$ together with a group homomorphism $Z\colon K(\D)\to\C$ having the property that
\[0\neq E\in\P(\phi)\implies Z(E)\in\R_{>0}\exp(i\pi\phi).\]
The map $Z$ is called the central charge, the nonzero objects of $\P(\phi)$ are called the semistables of phase $\phi$ and the simple
objects  are stable.
The smallest extension-closed subcategory of $\D$ containing the objects of $\P(\phi)$ for each $\phi\in(0,1]$ is an abelian category
called the heart of the stability
condition $\sigma$.

In fact $\sigma$ is completely determined by
its heart together with the central charge $Z$, and conversely, if $\AA\subset\D$ is the heart of a
bounded t-structure, and $Z\colon \K(\D)\to\C$
is a group homomorphism with the property that
\begin{equation*}
\tag{${*}$}
0\neq E\in\AA \implies Z(E)\in \R_{>0}\exp(i\pi\phi)\text{ with }\phi\in (0,1],
\end{equation*}
then there is a stability condition $\sigma$ on $\D$ with heart $\AA$ and central charge $Z$, providing that
$\AA$ satisfies a certain Harder-Narasimhan property with respect to $Z$. This property is automatically satisfied if $\AA$ has finite length.

Let us assume that $K(\D)$ is of finite rank. The set of all stability conditions satisfying a
technical condition called local-finiteness then form a complex manifold $\Stab(\D)$. There is a continuous forgetful map
\[\ZZ\colon \Stab(\D)\lra\Hom_{\Z}(K(\D),\C)\]
sending a stability condition to its central charge.

The group of exact autoequivalences $\Aut(\D)$ of $\D$ act on $\Stab(\D)$: an element $\Phi\in\Aut(\D)$ sends $(Z,\P)$ to $(Z',\P')$,
where $\P'(\phi)=\Phi(\P(\phi))$ and $Z'(E)=Z(\Phi^{-1}(E))$. The additive group $\C$ also acts on $\Stab(\D)$: an element $\lambda\in\C$ sends $(Z,\P)$
to $(Z',\P')$ where $\P'(\phi)=\P(\phi+\Re(\lambda))$ and $Z'(E)=\exp(-i\pi \lambda) Z(E)$. These two actions commute, and the action of the shift functor $[1]$ coincides with the action of $1\in\C$.

% *******************************************************************************************************************
% *******************************************************************************************************************
% *******************************************************************************************************************

\section{The results}

In this section we prove Theorem 1.3 and Corollary 1.4 as stated in the introduction. The corresponding results Theorem 1.1 and Corollary 1.2 for the Dynkin case
case are entirely analogous, but easier, and we shall confine ourselves to a few remarks on the proof at the end.

\subsection{}
Recall that we can identify the Grothendieck group $K(\Dh)$ with the affine root lattice $\Z\Lambdah\subset \hh^*$.
The classes $\alpha_i=[S_i]$ of the $n+1$ simple objects $S_i\in\AAh$ define a set of simple roots in $\Lambdah$.
The associated  Weyl chamber is the subset
\[\{\phi\in \Z\Lambdah\tensor \R: \phi(S_i)>0 \text{ for all }i\}.\]
Since $\Z\Lambdah\tensor \C=\hh^*$ we can identify group homomorphisms $Z\colon K(\Dh)\to \C$ 
with elements of $\hh$. Thus we have a continuous map
$\ZZ\colon \Stab(\Dh)\to \hh$.

\begin{lemma}
For each point $Z$ in the complexified Weyl chamber
\[R=\{Z\in \hh: \Im Z(\alpha_i)>0 \text{ for all }i\}\subset\hh\]
there is a unique stability condition $\sigma\in \Stab(\Dh)$ with heart $\AAh$ and central charge $Z$.
These points form a region $U\subset \Stab(\Dh)$ which is mapped homeomorphically by $\ZZ$ onto $R$.
\end{lemma}

\begin{pf}
The standard t-structure on $\D^b\Coh_G(\C^2)$ induces a bounded t-structure on $\Dh$ with heart $\AAh$.
Since $\AAh$ has finite length, the class of any nonzero element $E\in\AAh$ is a  positive
linear combination of the simple roots $\alpha_i$. Thus the condition $(*)$ holds. The Harder-Narasimhan property is automatic because $\AAh$ has finite length. The resulting stability condition
$\sigma=(Z,\P)$ is locally finite because for any $\phi$ there is an
$\epsilon>0$  such that the subcategories $\P((\phi-\epsilon,\phi+\epsilon))$ are contained in some shift of $\AAh$  and hence are of finite length.
\end{pf}

We shall need a result which follows from work of A. Craw and Ishii on moduli of $G$-constellations \cite[Proposition 2.2]{CI}. 
 
\begin{lemma}
\label{i}
For each point $\sigma\in U$ there is a semistable object whose class in $K(\Dh)=\Z\Lambdah$ is the imaginary root
$(0,1)$.
\end{lemma}

\begin{pf}
 A G-constellation is a representation of the group ring $\C[x,y]*G$ which as a $\C[G]$-module is isomorphic to $\C[G]$. Craw and Ishii use general
results of A.D. King \cite{Ki} to show that for any generic choice of weights the moduli space of semistable $G$-constellations is non-empty and has a projective morphism to the quotient $X=\C^2/G$. The fibre over the origin is then non-empty and consists of nilpotent representations. These define objects of $\AAh$ whose class in $K(\Dh)$ is the imaginary root $(0,1)$.
\end{pf}

Now we need two general results from \cite{Br2}. Let $\Stab_0(\Dh)$ be the connected component of $\Stab(\Dh)$ containing the subset $U$.
  
\begin{prop}
\label{cover}
The map $\ZZ:\Stab_0(\Dh)\to\hh$ sending a stability condition to its central
charge is a local homeomorphism onto an open subset of $\hh$ containing $\hregh$. The restriction
to $\ZZ^{-1}(\hregh)$ 
is a covering map.
\end{prop}

\begin{pf}
It is a general result \cite{Br1} that $\ZZ$ is a local homeomorphism onto an open subset of some linear subspace of $\hh$.
Since this subspace contains $R$ it must be the whole of $\hh$.
The fact that the restriction is a covering map is proved in exactly the same way as
the corresponding result on coherent sheaves on K3 surfaces \cite[Section 7]{Br2}.
This then implies that the image contains $\hregh$.
\end{pf}

\begin{prop}
Let $E\in\Dh$ be stable in a stability condition $\sigma\in\Stab_0(\Dh)$. Then  there is an open neighbourhood
$\sigma\in N\subset \Stab_0(\Dh)$ such that $E$ is stable for all stability conditions in $N$.
\end{prop}

\begin{pf}
Again this is exactly the same argument as in the K3 surface case, see \cite[Section 8]{Br2}.
In fact the only property of $\Stab_0(\Dh)$ needed is that it contains points
$\sigma$ satisfying $\ZZ(\sigma)\in\hregh$.
\end{pf}

\subsection{}

The following result shows that the autoequivalences $\Phi_{S_i}$ preserve the connected component $\Stab_0(\Dh)$
so that the group $\Br(\Dh)$ acts on $\Stab_0(\Dh)$.

\begin{lemma}
\label{l}
Let $\sigma=(Z,\P)$ be a point in the boundary of $U$ for which there is a unique
simple $S_i\in\AAh$ with $\Im Z(S_i)=0$. Assume further that $Z(S_i)\in \R_{<0}$.
Then the stability condition $\Phi_{S_i}^{-1}(\sigma)$ also lies in the boundary of $U$.
\end{lemma}

\begin{pf}
To help with notation set $T=S_i$. Take a small neighbourhood $V$ of $\sigma$ in $\Stab(\Dh)$ and consider the open subset
\[V_+=\{\sigma=(Z,\P)\in V: \Im Z(T)<0\}.\]
We claim that we can choose  $V$ small enough so that $\Phi_{T}^{-1} (V_+)\subset U$. It follows that the stability condition $\Phi_T^{-1}(\sigma)$ lies in the closure of $U$.
It cannot lie in $U$ because $Z(\Phi_T(T))=Z(T[-1])$ lies on the positive real axis.

Thus we are required to prove that if $V$ is small enough the heart of any $\tau\in V_+$ is equal to $\Phi_T(\AAh)\subset\Dh$.
It is a simple fact that if $\CC$ and $\CC'$ are both hearts of bounded t-structures in a triangulated category, and $\CC\subset\CC'$ then $\CC=\CC'$.
Since $\AAh$ has finite length it will therefore be enough to prove that  for all $j$ the object $\Phi_{T}(S_j)$ lies in the heart of any $\tau\in V_+$.

 Assume first that $j\neq i$. If vertices $j$ and $i$ are joined by an edge in $\Gammah$ then $\Hom^1_{\Dh}(S_i,S_j)=\C$ so there is a non-split short exact sequence
in $\AAh$ \[0\lra S_j\lra \Phi_T(S_j) \lra T\lra 0\]
It follows that $\Phi_T(S_j)$ is in the heart of $\sigma$ and its semistable factors have phases in the interval $(0,1)$. Choosing $V$ small enough, we can assume that this is the case for all $\tau\in V$ too. If $i$ and $j$ are not joined by an edge then $\Phi_T(S_j)=S_j$ and the same argument applies.

Finally consider $\Phi_T(T)=T[-1]$. Since $T$ was stable in $\sigma$ with phase 1, we can assume that $T$ is stable for all $\tau\in V$ too, with phase at most 2. Clearly one must have $\phi(T)>1$ for $\tau\in V_+$. This implies that $T[-1]$ has phase in the interval $(0,1]$ and hence lies in the heart of $\tau$. 
\end{pf}

Lemma \ref{l} shows that the autoequivalence $\Phi_{S_i}$ exchanges the two pieces of the boundary of $U$
given by $Z(S_i)\in \R_{<0}$ and $Z(S_i)\in\R_{>0}$. The crucial thing is that this autoequivalence reverses the orientations, taking the side where
$\Im Z(S_i)>0$ to the side where $\Im Z(S_i)<0$. This observation easily gives the following.

\begin{lemma}
\label{q}
For every stability condition $\sigma=(Z,\P)\in\Stab_0(\Dh)$ the central charge $Z$ does not vanish on the imaginary roots $(0,d)\in \Z\Lambdah\subset \hh^*$.  Furthermore, there is an
autoequivalence $\Phi\in\Br(\Dh)$ and an element $\lambda\in\C$ such that $\lambda \Phi(\sigma)$ lies in the closure of $U$.
\end{lemma}

\begin{pf}
First assume that the central charge $Z$ of $\sigma$ does not vanish on the imaginary roots $(0,d)\in\Z\Lambdah$.
Choose a path $\gamma$ joining $\sigma$ to a point of $U$. Since $\ZZ$ is a local homeomorphism we can assume that $Z((0,1))\neq 0$ for all stability conditions on the path $\gamma$. Normalising with the $\C$ action on $\Stab(\Dh)$ we can replace $\sigma$ by some $\lambda(\sigma)$ and assume that $\gamma$ lies in the affine slice
\[\hregh_{a}=\{(\theta,n)\in\hregh:n=i\}.\]
In this slice the complexified Weyl alcoves form a nice polyhedral decomposition, and since $\ZZ$ is a local homeomorphism we can wiggle the path $\gamma$ a bit so that it
passes through finitely many Weyl alcoves, and only passes through codimension one walls. Each time $\gamma$ passes through a wall Lemma \ref{l} shows that there is an element
of $\Br(\Dh)$ that takes one back to a stability condition in the closure of $U$. The result then follows.

Now suppose $Z((0,d))=0$. In particular, there are no semistable objects in $\sigma$ whose class in $K(\Dh)$ is the imaginary root $(0,1)$. By the results of \cite[Section 8]{Br2} this is true in an open neighbourhood of $\sigma$ in $\Stab_0(\Dh)$. But by the first part, and Lemma \ref{i}, stability conditions near $\sigma$ for which $Z((0,1))\neq 0$ do have semistable objects with class $(0,1)$. This gives a contradiction.
\end{pf}  

\subsection{}
Now we can prove our main results.

\medskip

{\noindent \bf Proof of Theorem 1.3.}
First we show that the image of the map $\ZZ$ is contained in $\hregh$.  Proposition \ref{cover} then shows that
$\ZZ$ is a covering map. Note that the autoequivalences $\Phi_{S_i}$ act on $K(\D)$ as root reflections in the simple roots $\alpha_i$. Thus the action of $\Br(\Dh)$ on $\Stab_0(\Dh)$ induces the action of the affine Weyl group on $\hh$ which preserves $\hregh$.
Similarly the action of $\C$ descends to the rescaling action of $\C^*$ which also preserves $\hregh$.

Suppose $\sigma=(Z,\P)\in \Stab_0(\Dh)$ satisfies $\ZZ(\sigma)\notin \hregh$. By the last result $Z((0,1))\neq 0$ so we can rescale and assume that $Z((0,1))=i$. Furthermore we can assume that $\sigma$ lies in the closure of $U$. Then the only roots which $Z$ can vanish on are the simple roots $\alpha_i$ defining the chamber $R$. These are the classes of the simple objects $S_i$ of $\AAh$. For all stability conditions in $U$ these objects $S_i$ are stable, so they are at least semistable in $\sigma$. It follows that the central charges $Z(S_i)$ do not vanish. 
 
By what was said above,
the  group $\Br(\Dh)$ acts as deck transformations for the covering map $\Stab_0(\Dh)\to\hregh/\Wh$.
Conversely, suppose two points $\sigma_1,\sigma_2\in \Stab_0(\Dh)$ map to the same point in $\hregh/\Wh$, and assume $\sigma_1\in U$.
Applying Lemma \ref{q} shows that there is an element $\Phi\in \Br(\Dh)$ and a $\lambda\in\C$ such that $\lambda\Phi(\sigma_2)\in \bar{U}$.
But since the complexified Weyl chamber $R$ is a fundamental domain for the action of $\Wh$ on $\hh$ which commutes with the $\C^*$-action it follows that $\lambda=2n$ is an even integer and $\sigma_1=\Phi(\sigma_2)[2n]$.
\qed
\medskip

{\noindent \bf Proof of Corollary 1.4.} 
Suppose an autoequivalence $\Psi$ of $\Dh$ preserves the connected component $\Stab_0(\Dh)$, take a stability condition $\sigma\in U$ and consider the stability condition $\Psi(\sigma)$. By Lemma \ref{q} there is an element $\Phi\in\Br(\Dh)$ and a $\tau \in \bar{U}$ such that
$\tau=\lambda \Phi \Psi(\sigma)$ for some $\lambda\in \C$.
Suppose we chose $\sigma\in U$ so that $\Re Z(S_i)=0$ for all $i$. Then there is some shift $[n]$ such that the stability condition $\sigma'=\lambda(\sigma)[-n]$ lies in $U$. Now $\Upsilon=\Phi\Psi[n]$ takes $\sigma'$ to $\tau$. Deforming $\sigma'$ and $\tau$ a little bit we can assume that they both lie in $U$ and hence have heart $\AAh$. It follows that $\Upsilon$ fixes $\AAh\subset \Dh$.

Now $\Upsilon$ permutes the simple objects of $\AAh$ and hence induces an automorphism of the graph $\Gammah$. Viewing $\Dh$ as a subcategory of the derived category of representations of the preprojective algebra of $\Gammah$ it is easy to see that conversely any automorphism of $\Gammah$ lifts to an automorphism of $\Dh$ preserving $\AAh$. Thus we may assume that $\Upsilon$ fixes the simples $S_i$. But then it acts trivially on $K(\Dh)$ and hence fixes all stability conditions $\sigma\in U$. It follows that it acts trivially on $\Stab_0(\Dh)$.
\qed

\medskip
The proofs in the finite type cases proceed in exactly the same way. Each result we proved in this section holds also for the Dynkin case on doffing hats and replacing the
$n+1$ simples $S_0,\cdots,S_n$
of $\AAh$ with the $n$ simples $S_1,\cdots,S_n$ of $\AA$. The only exception is
Lemma \ref{i} which is meaningless in the Dynkin case. The proof of Lemma \ref{q} is easier since the decomposition of $\h$ into Weyl chambers is polyhedral so we have
no need to rescale or pass
to an affine slice.

\begin{remark}\label{eend}
In 2019 Michael Wemyss pointed out a gap in the published version of this paper. Namely, nowhere is it proved that the subgroup $\Br(\Dh)$ acts faithfully on the space $\Stab_0(\Dh)$, nor that $\Br(\D)$ acts faithfully on $\Stab_0(\D)$. We now  briefly sketch an argument for these two statements, beginning with the affine case. 

Note first that the  twist functors $\Phi_{S_i}$ can be realised as auto-equivalences of $\D^b\Coh(Y)$ given by Fourier-Mukai functors, and hence the same is true of all elements of the subgroup $\Br(\Dh)$. Moreover,  these auto-equivalences fix all skyscraper sheaves of points lying off the exceptional locus of the contraction $f\colon Y\to X$.
 
 Suppose that some element $\Psi\in \Br(\Dh)$ acts trivially on $\Stab_0(\Dh)$. Then $\Psi$  fixes the standard heart $\AAh\subset\Dh$ and acts trivially on the Grothendieck group $K(\Dh)$, and hence fixes the simple objects of $\AAh$ pointwise. It follows that it permutes the skyscraper sheaves of points on the exceptional locus of $f$, since for a suitable stability condition on $\AAh$, these are the stable objects of the given class. Thus $\Psi$ is a Fourier-Mukai transform which takes all skyscrapers to skyscrapers.
 
 Any such auto-equivalence $\Psi$ is necessarily of the form $\Psi(-)=g^*(-)\tensor L$ for some automorphism $g\colon Y\to Y$, and some line bundle $L\in \operatorname{Pic}(Y)$. It is clear that $g$ must be the identity, since it is the identity on the complement of the exceptional locus of $f$. Since $\Psi$ preserves the simple objects of $\AAh$, the line bundle $L$ is trivial when restricted to each irreducible component of the exceptional locus. It follows that it is trivial, and so $\Psi$ is the identity.
 
 Let us now consider  the finite-type case. Suppose some element $\Psi\in \Br(\Dh)$ acts trivially on $\Stab_0(\Dh)$. Then as before, $\Psi$ fixes the simple objects of the heart $\AA$ pointwise. Note that as a composition of the twist functors $\Phi_{S_i}$, we can view the equivalence $\Psi$ as the restriction of an element $\Psi\in \Br(\Dh)$, which commutes with the  functor $\R f_*\colon \D^b\Coh(Y)\to \D^b\Coh(X)$. After what was proved above it remains to show that $\Psi$ fixes the simple $S_0\in \AAh$ corresponding to the extending vertex of the affine Dynkin diagram. 
 
Set $E=\Psi(S_0)\in \Dh$. Note that the functor $\R f_*$ is exact with respect to the t-structures whose  hearts  are $\AAh$ and $\Coh(Y)$. Since $\R f_*(E)=\R f_*(S_0)\in \Coh(X)$ it follows that for each $i\neq 0$ the  cohomology object $H^i(E)\in \AAh$ satisfies $\R f_*(H^i(E))=0$, and hence $H^i(E)\in \AA$.  Now $\Psi$ preserves $\AA$, so a map $A[i]\to E$ with $A\in \AA$ and $i>0$ gives a map $A'[i]\to S_0$ with $A'=\Psi^{-1}(A)\in \AA$, which is impossible. This shows that $H^i(E)=0$ for $i<0$, and a similar argument shows that $H^i(E)=0$ for $i>0$. We therefore conclude that $E\in \AAh$, and hence that $\Psi$ preserves the heart $\AAh$. Since $\AAh$ is finite length and $\Psi$ is an equivalence it is then easy to see that $\Psi (S_0)=S_0$, which completes the argument.
\end{remark}

\bigskip

\noindent Department of Pure Mathematics,
University of Sheffield,
Hicks Building, Hounsfield Road, Sheffield, S3 7RH, UK.

\smallskip

\noindent email: {\tt t.bridgeland@sheffield.ac.uk}

\end{document}